\documentclass{article}

\usepackage{amssymb,amscd,amsmath}
\usepackage[matrix, arrow]{xy}

\newtheorem{Thm}{Theorem}
\newtheorem{Lem}[Thm]{Lemma}
\newtheorem{Prop}[Thm]{Proposition}
\newtheorem{Cor}[Thm]{Corollary}

\newenvironment{Pf}{\textbf{Proof }}{\hspace{\stretch{1}}$\square$}

\DeclareMathOperator{\dimv}{\underline\dim}
\DeclareMathOperator{\Hom}{Hom}
\DeclareMathOperator{\Ext}{Ext}
\DeclareMathOperator{\End}{End}
\DeclareMathOperator{\Aut}{Aut}

\DeclareMathOperator{\HA}{\mathcal H}
\DeclareMathOperator{\CA}{\mathcal C}

\DeclareMathOperator{\QG}{\mathcal U}

\newcommand{\bQ}{\mathbb Q}

\newcommand{\diff}{\frac{\mathrm d}{\mathrm dt}}

\begin{document}

\begin{center}
{\Large\bf The Composition Algebra of an Affine Quiver}
\vspace{\baselineskip}

Andrew Hubery\\
Universit\"at Paderborn, Germany\\
hubery@math.uni-paderborn.de
\end{center}

\begin{abstract}
We study the Hall and composition algebras of an affine quiver. In the case of a cyclic quiver, we provide generators for the central polynomial algebra described by Schiffmann and prove that this is in fact the whole of the centre of the Hall algebra. For an affine quiver without oriented cycles, we obtain a structure theorem for the composition algebra refining the structure provided by Zhang.
\end{abstract}


\section{Introduction}

There are many interesting connections between the representation theory of quivers and the structure theory of Kac-Moody Lie algebras, with one of the main tools being that of the Ringel-Hall algebra \cite{Ringel1}, or equivalently the quantised enveloping algebra \cite{Lusztig}.

In the case of a Dynkin quiver, this is particularly well understood (see for example \cite{Lusztig1, Reineke, Ringel3}). For an affine quiver, though, despite the explicit description of the module category provided by \cite{DR}, a concrete description of the structure of the composition algebra remained open.

Here we provide such a structure. We obtain a refinement of the triangular decomposition $\mathcal C=\mathcal P\otimes\mathcal T\otimes\mathcal I$ given by Zhang \cite{Zhang}, corresponding to the preprojective, regular and preinjective modules respectively. In particular, we show in Theorem \ref{ThmC} that the subalgebra $\mathcal T$ can be written as the product of an infinitely generated polynomial ring together with the composition algebras of the non-homogeneous tubes. We remark that this also generalises the structure for the Kronecker quiver provided by Sz\'ant\'o \cite{Szanto}. Note, however, that our proof is entirely independent.

We remark that this result reduces the problem of describing the elements in the composition algebra of an arbitrary affine quiver to that of the composition algebra of a non-homogeneous tube, or equivalently the composition algebra of nilpotent representations of a cyclic quiver. This has been studied by Ringel \cite{Ringel2} and more recently by Deng and Du \cite{DD}.

Regarding the Hall algebra of a cyclic quiver is also considered, we know from the work of Schiffmann \cite{Schiff} that there is a central polynomial subalgebra on infinitely many generators which, together with the composition algebra, generates the Hall algebra. We describe a set of generators for this polynomial ring in terms of nilpotent representations and prove that they in fact generate the whole of the centre. In particular, the composition algebra of a cyclic quiver has trivial centre.

The paper is organised as follows. We first recall some general results about Hall algebras and prove that the composition algebra for an arbitrary quiver without oriented cycles contains all the preprojective and preinjective modules.

We then consider the classical isomorphism between the Hall algebra of nilpotent representations over the quiver with one vertex and one loop and Macdonald's ring of symmetric functions. We list some important generators for this algebra and describe their generating functions.

In Section 4 we study an arbitrary cyclic quiver. Here we prove our first main result, Theorem \ref{ThmA}, concerning the generators for Schiffmann's central subalgebra of the Hall algebra, and that this is in fact the whole of the centre.

Finally in Section 5 we study the affine quivers without oriented cycles. We first obtain in Theorem \ref{ThmB} certain regular elements contained in the composition algebra. This result is a generalisation of Theorem 4.3 in \cite{Szanto}. Using this, we then prove Theorem \ref{ThmC}, which describes the structure of the composition algebra and provides us with a PBW basis.

\textit{Acknowledgements.} The author would like to thank Prof. C.M.~Ringel for his interest and support.

\section{General results about Hall algebras}

Let $Q$ be a quiver without oriented cycles and let $k$ be a finite field.

The Hall algebra $\HA(kQ)$ \cite{Ringel1}is the $\bQ$-algebra with basis the isomorphism classes $[M]$ of finite dimensional modules for $kQ$ and multiplication
$$[M][N]:=\sum_{[X]}F_{MN}^X[X], \quad\textrm{where}\quad F_{MN}^X:=\left|\{Y\subset X:X/Y\cong M, Y\cong N\}\right|.$$
The algebra $\HA(kQ)$ is naturally graded by dimension vector.

The composition algebra $\CA(kQ)$ is defined to be the subalgebra generated by the simple modules. It is clearly a graded subalgebra.

Consider now the quantised enveloping algebra $\QG_q^+(\mathfrak g)$ of the Kac-Moody Lie algebra with generalised Cartan matrix that of the quiver $Q$ \cite{Lusztig}. Twisting the multiplication by the Euler form of $Q$, we obtain an algebra $\mathcal C'_q(Q)$. We note that this is a $\bQ(q)$-algebra generated by elements $e_i$ for $i$ a vertex of $Q$ and with defining relations the quantum Serre relations. For example, if the vertices $i$ and $j$ are not connected by an arrow, then the generators $e_i$ and $e_j$ commute. If there is a single arrow $i\to j$, then we have the relations
$$e_i^2e_j-(q+1)e_ie_je_i+qe_je_i^2 \quad\textrm{and}\quad e_ie_j^2-(q+1)e_je_ie_j+qe_j^2e_i.$$
Now let $\mathcal C_q(Q)$ be the $\bQ[q]$-subalgebra of $\mathcal C'_q(Q)$ generated by the vertices.

From the work of Green \cite{Green}, Lusztig \cite{Lusztig} and Ringel \cite{Ringel1} (see also \cite{Ringel4}) we know that the specialisation of $\CA_q(Q)$ at $q=\left|k\right|$ is isomorphic to the composition algebra $\CA(kQ)$, whereas the specialisation to $q=1$ recovers the universal enveloping algebra $\QG^+(\mathfrak g)$ of the positive part of $\mathfrak g$. For this reason, $\CA_q(Q)$ is called the generic composition algebra of $Q$.

We also deduce that the dimensions $\dim\CA_q(Q)_\alpha$ and $\dim\QG_q^+(\mathfrak g)_\alpha$ coincide, and both equal $\dim\QG^+(\mathfrak g)_\alpha$. In particular, this dimension can be determined using a PBW basis of $\QG^+(\mathfrak g)$.

\subsection{Reflection functors}\label{RF}

Let $i$ be a sink of $Q$ and let $\widetilde Q:=r_iQ$ be the quiver obtained by reversing all arrows involving $i$.

The reflection functors
$$R_i^+:(\bmod kQ)\langle i\rangle\to(\bmod k\widetilde Q)\langle i\rangle \quad\textrm{and}\quad R_i^-:(\bmod k\widetilde Q)\langle i\rangle\to(\bmod kQ)\langle i\rangle$$
are inverse equivalences between the subcategories of all modules not containing the simple $S_i$ as a direct summand \cite{DR}.

These functors naturally induce algebra isomorphisms
$$R_i^+:\HA(kQ)\langle i\rangle\to\HA(k\widetilde Q)\langle i\rangle \quad\textrm{and}\quad R_i^-:\HA(k\widetilde Q)\langle i\rangle\to\HA(Q)\langle i\rangle.$$

\begin{Prop}
The subalgebra $\CA(kQ)\langle i\rangle$ is generated by the elements
$$\sigma(j)_a:=\sum_{\substack{\dimv X=(1,a)\\X\textrm{ indec}}}[X] \quad\textrm{for }0\leq a\leq n_j$$
for each subquiver of the form $j\xrightarrow{(n_j)}i$ ($n_j\geq 0$).
\end{Prop}

\begin{Pf}
Suppose that $Q$ is the quiver $j\xrightarrow{(n)}i$ for some $n\geq 0$. Then for each $r\geq 0$ we have that
$$\smash[t]{[S_j][S_i^r]=\sum_{a=0}^r[S_i^{r-a}]\sigma(j)_a.}$$
It follows by induction that each $\sigma(j)_a$ lies in $\CA(kQ)$. We note that, in particular, $\sigma(j)_0=[S_j]$ and $\sigma(j)_n=[P_j]$, where $P_j$ is the projective cover of $S_j$.

It follows that we have the decomposition
$$\smash[b]{\CA(kQ)=\sum_{r\geq 0}[S_i]^r\bQ\langle\sigma_0,\ldots,\sigma_n\rangle,}$$
and in particular,
$$\CA(kQ)\langle i\rangle = \bQ\langle\sigma_0,\ldots,\sigma_n\rangle.$$

For a general quiver $Q$, we deduce that we have the decomposition
$$\CA(kQ)=\sum_{r\geq 0}[S_i]^r\CA(kQ)\langle i\rangle,$$
with $\CA(kQ)\langle i\rangle$ generated by the elements $\sigma(j)_a$ for $0\leq a\leq n_j$ where $\smash{j\xrightarrow{(n_j)}i}$ (including the case $n_j=0$).
\end{Pf}

\begin{Prop}
The reflection functors induce isomorphisms between the algebras $\CA(kQ)\langle i\rangle$ and $\CA(k\widetilde Q)\langle i\rangle$.
\end{Prop}

\begin{Pf}
For each subquiver of the form $j\xrightarrow{(n)}i$, the functor $R_i^+$ takes $\sigma_a$ to $\tilde\sigma_{n-a}$, where
$$\tilde\sigma_a:=\sum_{\substack{\dimv X=(1,a)\\X\textrm{ indec}}}[X]$$
is defined analogously for the quiver $j\xleftarrow{(n)}i$.
\end{Pf}

\begin{Thm}\label{preproj}
The composition algebra contains the subalgebra $\mathcal P(kQ)$ (respectively $\mathcal I(kQ)$) generated by the isomorphism classes of all preprojective (respectively preinjective) modules.
\end{Thm}

\begin{Pf}
It is enough to show that the isomorphism class of each indecomposable preprojective lies in $\CA(kQ)$.

Let $P$ be an indecomposable preprojective. Then (by considering the Coxeter transformation) there exists a sequence of reflection functors such that $P':=R_{i_n}^+\cdots R_{i_i}^+(P)$ is simple for the quiver $r_{i_n}\cdots r_{i_1}Q$. Hence $[P']$ lies in the composition algebra for this quiver. By applying the functor $R_{i_1}^-\cdots R_{i_n}^-$, we deduce that $[P]\in\CA(kQ)$.
\end{Pf}

\section{The classical Hall algebra}\label{HA}

The reference for this section will be Macdonald's book \cite{Macd}.

We fix a finite field $k$ and let $d\geq1$. Let $K/k$ be a field extension of degree $d$.

Let $\HA_d=\HA_d(\widetilde{\mathbb A}_{\,0})$ be the generic Hall algebra (over $\bQ(q)$) arising from the Hall algebra of nilpotent representations for the path algebra $K\widetilde{\mathbb A}_{\,0}$. That is, $\HA_d$ is the $\bQ(q)$-algebra with basis $u_d(\lambda)$ indexed by partitions $\lambda$ and with multiplication
$$u_d(\lambda)u_d(\mu)=\sum_\nu F_{\lambda\mu}^\nu(q^d)u_d(\nu),$$
where the $F_{\lambda\mu}^\nu(q)$ are the classical Hall polynomials.

Let $\Gamma=\bQ(q)[X_1,X_2,\ldots]^{\mathfrak S_\infty}$ be Macdonald's ring of symmetric functions. This contains the Hall-Littlewood polynomials $P_\lambda(X;s)$, which are polynomials in $\Gamma[s]$ indexed by partitions $\lambda$. For a partition $\lambda=(\lambda_1\geq\lambda_2\geq\cdots)$ we define $n(\lambda):=\sum_i(i-1)\lambda_i$ and set $\ell(\lambda)$ to be the length of $\lambda$, i.e. the number of non-zero parts $\lambda_i$.

The following theorem was first conjectured by Steineitz \cite{Stein}, but is generally attributed to Hall \cite{Hall}. See also \cite{Macd}.
\begin{Thm}
The map $\psi_d:\Gamma\to\HA_d$ sending $P_\lambda(X;q^{-d})$ to $q^{n(\lambda)d}u_d(\lambda)$ is an algebra isomorphism.
\end{Thm}

We now consider various kinds of elements of $\Gamma$ and their images under $\psi_d$.

The elementary symmetric functions $e(n)=P_{(1^n)}(X;s)$ are given by the generating function
$$E(t)=1+\sum_{n\geq1}e(n)t^n=\prod_i(1+X_it).$$
Applying $\psi_d$ we obtain
$$\psi_d(e(n))=q^{dn(n-1)/2}u_d(1^n)$$
and we note that $u_d(1^n)$ is the isomorphism class of the elementary module of dimension $n$.

The complete symmetric functions $h(n)$ are given by the generating function
$$H(t)=1+\sum_{n\geq1}h(n)t^n=\prod_i(1-X_it)^{-1}.$$
Applying $\psi_d$ and using III.3 Example 1 of \cite{Macd} we obtain
$$\psi_d(h(n))=\sum_{\lambda\vdash n}u_d(\lambda).$$

The power sum functions $p(n)$ are given by the generating function
$$P(t)=\sum_{n\geq1}p(n)t^{n-1}=\diff\log H(t).$$
Using III.7 Example 2 of \cite{Macd} we obtain
$$\psi_d(p(n))=\sum_{\lambda\vdash n}(1-q^d)\cdots(1-q^{(\ell(\lambda)-1)d})u_d(\lambda).$$

Now set $c(s,n):=(1-s)P_{(n)}(X;s)$ (denoted $q_n(X;s)$ in \cite{Macd}). These have the generating function
$$C(s,t)=1+\sum_{n\geq1}c(s,n)t^n=H(t)/H(st).$$
In particular,
$$\diff\log C(s,t)=\diff\big(\log H(t)-\log H(st)\big)=P(t)-sP(st)=\sum_{n\geq1}(1-s^n)p(n)t^{n-1}.$$
Thus
$$\psi_d(c(q^{-d},n))=(1-q^{-d})u_d(n)$$
and $u_d(n)$ is the isomorphism class of the cyclic module of dimension $n$ (hence the notation $c(s,n)$).

We note that $\Gamma$ is an infinite polynomial ring over $\bQ(q)$ with respect to any of the sets of generators $e(n)$, $h(n)$ or $p(n)$. Moreover, it follows from the results in III.2 in \cite{Macd} that the same holds for the $c(n,s)$ whenever $s\in\bQ(q)$ is not a root of unity (i.e. $s\neq\pm 1$).

We shall need the following result. The automorphism group of the cyclic module of dimension $n$ has size
$$a_d(n)=(q^d-1)q^{(n-1)d}=(1-q^{-d})q^{nd}.$$
The generating function $A_d(t):=1+\sum_{n\geq1}a_d(n)u_d(n)t^n$ therefore satisfies
$$A_d(t)=1+\sum_{n\geq1}(1-q^{-d})u_d(n)(q^dt)^n=\psi_d(C(q^{-d},q^dt)).$$

In particular,
\begin{align*}
\diff\log A_d(t) &= \psi_d\Big(\diff\log C(q^{-d},q^dt)\Big)=\psi_d\Big(q^d\sum_{n\geq1}(1-q^{-dn})p(n)(q^dt)^{n-1}\Big)\\
&= \sum_{n\geq1}(q^{nd}-1)\psi_d(p(n))t^{n-1}.
\end{align*}

\section{The cyclic quiver case}

We shall now consider the Hall algebra of nilpotent representations of the cyclic quiver $\widetilde{\mathbb A}_{\,l}$. This has vertices labelled $0$ to $l$ (with the convention that $l+1=0$) and arrows $i\to i+1$.

The generic composition algebra $\CA_q(\widetilde{\mathbb A}_{\,l})$ is again related to the quantised enveloping algebra $\QG_q^+(\widehat{\mathfrak{sl}}_{\,l+1})$ by twisting the multiplication using the Euler form \cite{Ringel2}.

We recall a result of Schiffmann \cite{Schiff}.
\begin{Thm}\label{SchiffThm}
The generic Hall algebra of nilpotent representations of the cyclic quiver of type $\widetilde{\mathbb A}_{\,l}$ can be decomposed as
$$\HA_q(\widetilde{\mathbb A}_{\,l})=\CA_q(\widetilde{\mathbb A}_{\,l})\otimes\bQ(q)[z_1,z_2,\ldots],$$
where $\CA(\tilde{\mathbb A}_{\,l})$ is the composition algebra and $z_m$ is central of degree $m\delta$.
\end{Thm}

Let $\Delta$ denote Green's comultiplication and $\{-,-\}$ the symmetric bilinear form \cite{Green}. Therefore
$$\Delta([X])=\sum_{[M],[N]}F_{MN}^X\frac{\left|\Aut M\right|\left|\Aut N\right|}{\left|\Aut X\right|}[M]\otimes[N]$$
and
$$\{[M],[N]\}=\delta_{[M][N]}\frac{1}{\left|\Aut M\right|}.$$
These are related via $\{x,yz\}=\{\Delta(x),y\otimes z\}$.

The central subalgebra $\bQ(q)[z_1,z_2,\ldots]$ was characterised by Schiffmann as the intersection of the $\mathrm{Ker}(e_i')$ for all vertices $i$, where the operator $e_i'$ is defined to be adjoint with respect to $\{-,-\}$ to premultiplication by the simple $[S_i]$. That is,
$$\{e_i'x,y\}=\{x,[S_i]y\}=\{\Delta(x),[S_i]\otimes y\}.$$

The main result of this section is the following theorem.

\begin{Thm}\label{ThmA}
We can take as generators
$$z_m=\sum_{\substack{[M]:\dimv M=m\delta\\\mathrm{soc}M\textrm{ square free}}}(-1)^{\dim\End(M)}\left|\Aut M\right|[M].$$
Moreover, these generate the whole of the centre of $\HA(\widetilde{\mathbb A}_{\,l})$. In particular, the centre of $\CA(\widetilde{\mathbb A}_{\,l})$ is trivial.
\end{Thm}

We first note that the Auslander-Reiten translate $\tau$ satisfies ${}^\tau S_i=S_{i+1}$. Using the Euler form, we then see that
\begin{align*}
\langle M,S_i\rangle &= \dim\Hom(M,S_i)-\dim\Ext^1(M,S_i)\\
&= \dim\Hom(M,S_i)-\dim\Hom(S_{i-1},M).
\end{align*}
In the special case that $\dimv M=m\delta$, $\langle M,-\rangle$ is identically zero and thus the socle is square free if and only if the top is square free.

We divide the proof of the theorem into several parts.

\begin{Prop}
Each $z_m$ lies in $\mathrm{Ker}(e_i')$ for every vertex $i$.
\end{Prop}

\begin{Pf}
We wish to show that no term of the form $[S_i]\otimes[N]$ occurs in $\Delta(z_m)$.

Suppose there exists a short exact sequence of the form
$$0\to N\to X\to S_i\to 0$$
with $\dimv X=m\delta$ and $\mathrm{soc}X$ square free. Then we must have
$$\dim\Hom(N,S_i)=0, \quad \dim\Hom(N,S_{i+1})\leq 2 \quad\textrm{and}\quad F_{S_iN}^X=1.$$

Let us write $N=S_{i+1}(r)\oplus S_{i+1}(s)\oplus N'$ with $\dim\Hom(N',S_i\oplus S_{i+1})=0$, where $S_i(a)$ is the unique indecomposable of length $a$ and top $S_i$. Since $\mathrm{soc}N$ is square free, $r\not\equiv s\mod l+1$ and so we may assume that $r>s$.

Set $M_1:=S_i(r+1)\oplus S_{i+1}(s)$ and $M_2:=S_{i+1}(r)\oplus S_i(s+1)$. Then clearly $M_1\oplus N'$ and $M_2\oplus N'$ are the only possible choices for $X$. Thus it is enough to prove that the dimensions of the endomorphism rings of these two modules differ by 1.

We note the following equalities.
\begin{enumerate}
\item $\dim\Hom(L,S_i(a+1))=\dim\Hom(L,S_{i+1}(a))$ for all $a\geq 0$ and all modules $L$ such that $\Hom(L,S_i)=0$.
\item $\dim\Hom(S_i(a+1),L)=\dim\Hom(S_{i+1}(a),L)+\dim\Hom(S_i,L)$ for all $a\geq0$ and all modules $L$ such that $\Hom(L,S_{i+1})=0$;
\item $\dim\Hom(S_{i+1}(a),S_i(b+1))=\dim\Hom(S_{i+1}(a),S_{i+1}(b))$ for all $a$ and $b$;
\item $\dim\Hom(S_i(a),S_{i+1}(b))=\begin{cases}\dim\Hom(S_i(a),S_i(b+1)) & \textrm{if } a\leq b;\\ \dim\Hom(S_i(a),S_i(b+1))-1 & \textrm{if } a>b\end{cases}$;
\end{enumerate}

Since $\Hom(N',S_i\oplus S_{i+1})=0$, we see that
$$\dim\Hom(N',M_1)=\dim\Hom(N',S_{i+1}(r)\oplus S_{i+1}(s))=\dim\Hom(N',M_2).$$
and that
\begin{align*}
\dim\Hom(M_1,N') &= \dim\Hom(S_{i+1}(r)\oplus S_{i+1}(s),N')+\dim\Hom(S_i,N')\\
&= \dim\Hom(M_2,N').
\end{align*}
Similarly, we deduce that $\dim\End(M_1)$ equals
$$\dim\End(S_i(r+1)\oplus S_i(s+1))-\dim\Hom(S_i,S_i(r+1)\oplus S_i(s+1))-1$$
whereas $\dim\End(M_2)$ equals
$$\dim\End(S_i(r+1)\oplus S_i(s+1))-\dim\Hom(S_i,S_i(r+1)\oplus S_i(s+1)),$$
and so we are done.
\end{Pf}

Given a module $M$, we can consider the Loewy lengths of its indecomposable summands. These will then determine a partition, denoted $\mu(M)$.

For partitions $\lambda$ and $\mu$ we define their cup product $\lambda\cup\mu$ to be the partition formed by arranging all the parts of $\lambda$ and $\mu$ in descending order. Moreover, if $\lambda$ and $\mu$ are both partitions of $n$ for some $n$, then we set $\lambda<\mu$ if the first time that $\lambda_i\neq\mu_i$ implies that $\lambda_i>\mu_i$ (i.e. the reverse lexicographic ordering).

\begin{Lem}\label{ext}
Any extension $X$ of $M$ by $N$ must satisfy $\mu(X)\leq\mu(M)\cup\mu(N)$. Moreover, we have equality if and only if $X\cong M\oplus N$.
\end{Lem}

\begin{Pf}
We shall prove the result by induction on the number of indecomposable summands of $M$.

Suppose that $M$ is indecomposable and consider a short exact sequence
$$\xymatrix{0\ar[r] &N\ar[r]^f &X\ar[r]^g &M\ar[r] &0.}$$
Let $X_i$ be the indecomposable summands of $X$ and write $g_i:X_i\to M$ for the restriction of $g$. Let $M_i$ be the image of $g_i$. By applying a suitable automorphism, we may assume that $X_i=\begin{smallmatrix}M_i\\A_i\end{smallmatrix}$ and that $g_i$ factors through the canonical maps $X_i\twoheadrightarrow M_i\hookrightarrow M$.

Suppose that there exists a morphism $\pi:X_j\to X_i$ such that $g_j=g_i\pi$. Then we can use the automorphism $\left(\begin{smallmatrix}1&\pi\\&1\end{smallmatrix}\right)$ of $X_i\oplus X_j$ to reduce to the case $g_j=0$. In particular, we may assume that the non-zero $M_i$ are all distinct.

Hence we can write $X=X_1\oplus\cdots\oplus X_r\oplus X'$ and $N=N'\oplus X'$ with $0\to N'\to X_1\oplus\cdots\oplus X_r\to M\to 0$ exact and $M=M_1>\cdots>M_r>0$. Also, we may assume that no $g_j$ factors through another $g_i$.

Now clearly each $A_i=\mathrm{Ker} g_i$ is indecomposable and they all have isomorphic tops. Therefore for $i<j$ either $A_i\twoheadrightarrow A_j$ or $A_j\twoheadrightarrow A_i$. In the latter case, we deduce that $g_j$ factors as $X_j\to X_i\to M$, a contradiction. Hence we have proper epimorphisms $A_i\twoheadrightarrow A_j$ for all $i<j$.

It follows that the kernel $N'$ equals $N_1\oplus\cdots\oplus N_r$, where $N_i=\begin{smallmatrix}M_{i+1}\\A_i\end{smallmatrix}$ (setting $M_{r+1}=0$). For, we have the morphism $\left(\begin{smallmatrix}f_i\\-f_i'\end{smallmatrix}\right):N_i\to X_i\oplus X_{i+1}$ where $f_i$ is the canonical monomorphism and $f_i'$ the canonical epimorphism. This determines a monomorphism $f:N'\to X_1\oplus\cdots\oplus X_r$ such that $gf=0$, and clearly $\dimv N=\dimv X-\dimv M$.

Finally, let $m_i$ (respectively $a_i$) denote the Loewy length of $M_i$ (respectively $A_i$). Then
$$\mu(X)=(m_1+a_1,\ldots,m_r+a_r)\cup\mu(X')$$
whereas
$$\mu(N)=(m_2+a_1,\cdots,m_r+a_{r-1},a_r)\cup\mu(X') \quad\textrm{and}\quad \mu(M)=(m_1).$$
Suppose that $a_1=0$. Then $r=1$, $\mu(X)=\mu(M)\cup\mu(N)$ and $X\cong M\oplus N$. On the other hand, if $a_1\neq0$, then $\mu(X)<\mu(M)\cup\mu(N)$.

This proves the result when $M$ is indecomposable.

Suppose now that $M=M'\oplus M''$ with $M''$ indecomposable. Consider the exact commutative diagram
$$\xymatrix{0\ar[r] &N\ar[r]\ar@{=}[d] &Y\ar[r]\ar[d] &M'\ar[r]\ar[d] &0\\
0\ar[r] &N\ar[r] &X\ar[r]\ar[d] &M\ar[r]\ar[d] &0\\
&&M''\ar@{=}[r] &M''}$$
By induction we have
$$\mu(X)\leq\mu(Y)\cup\mu(M'') \quad\textrm{and}\quad \mu(Y)\leq\mu(M')\cup\mu(N),$$
and clearly $\mu(M)=\mu(M')\cup\mu(M'')$. Since the cup product is associative, the result follows.
\end{Pf}

\begin{Prop}
The $z_m$ are algebraically independent.
\end{Prop}

\begin{Pf}
Consider an arbitrary algebraic equation $\sum_{\lambda\vdash n}c_\lambda z_\lambda$, where $z_\lambda=\prod_t z_{\lambda_t}$.

Each module $M$ occurring in the expression for $z_m$ satisfies $\mu(M)\leq m^{\cup(l+1)}$, the partition consisting of $l+1$ copies of $m$. Moreover, there is a unique such module $Z_m$ such that $\mu(Z_m)=m^{\cup(l+1)}$, namely $Z_m=S_0(m)\oplus\cdots\oplus S_l(m)$.

It follows from Lemma \ref{ext} that every module $M$ occurring in the product $z_\lambda$ must satisfy $\mu(M)\leq\lambda^{\cup(l+1)}$, the partition given by repeating each term $\lambda_t$ of $\lambda$ $l+1$ times. Moreover, there is a unique summand giving rise to the partition $\lambda^{\cup(l+1)}$, namely $Z_\lambda=\bigoplus_t Z_{\lambda_t}$. Proceeding by induction on $\lambda$ we deduce that each $c_\lambda=0$ and hence that the $z_m$ are algebraically independent.
\end{Pf}

\begin{Lem}
For a module $M$ we have that $\dim\Hom(M,N)=\dim\Hom(N,M)$ for all $N$ if and only if $M\cong Z_\pi$ for some partition $\pi$. In particular, $\dimv M$ is a multiple of $\delta$ and $M\cong{}^\tau M$.
\end{Lem}

\begin{Pf}
We note that
$$\langle M,N\rangle=\dim\Hom(M,N)-\dim\Hom(N,{}^\tau M).$$
If $M\cong Z_\pi$, then $M\cong{}^\tau M$, $\dimv M$ is a multiple of $\delta$ and $\langle M,-\rangle$ is identically zero. Thus the condition is clearly sufficient.

Let $m_i(r)$ denote the multiplicity of $S_i(r)$ in $M$.

Consider the short exact sequence $0\to S_t\to S_0(t+1)\to S_0(t)\to 0$.

For an indecomposable $I$, we can lift a homorphism $I\to S_0(t)$ to $I\to S_0(t+1)$ unless $I\to S_0(t)$ is a monomorphism, and hence $I\cong S_j(t-j)$ for some $0\leq j<t$. It follows that $\dim\Hom(I,S_0(t+1))$ equals
$$\dim\Hom(I,S_t)+\dim\Hom(I,S_0(t))-\begin{cases}1 &\textrm{if }I\cong S_j(t-j), \textrm{ some } 0\leq j<t;\\
0 &\textrm{otherwise}.\end{cases}$$
We deduce that $\dim\Hom(M,S_0(t+1))$ equals
$$\dim\Hom(M,S_t)+\dim\Hom(M,S_0(t))-\big(m_0(t)+\cdots+m_{t-1}(1)\big).$$
Similarly, we have that $\dim\Hom(S_0(t+1),M)$ equals
$$\dim\Hom(S_t,M)+\dim\Hom(S_0(t),M)-\big(m_1(t)+\cdots+m_t(1)\big).$$

Suppose that $\dim\Hom(M,N)=\dim\Hom(N,M)$ for all $N$. Then, by the above when $t=1$, we deduce that $m_0(1)=m_1(1)$. We conclude that $m_i(1)=m(1)$ is constant for all vertices $i$. Now, by induction, we see that $m_i(t)=m(t)$ is constant for all vertices $i$, and for all $t\geq1$.

Hence there exists a partition $\pi$ such that $M\cong Z_\pi$.
\end{Pf}

\begin{Prop}
The centre of $\HA(\widetilde{\mathbb A}_{\,l})$ is precisely the subalgebra $\bQ(q)[z_1,\ldots]$.
\end{Prop}

\begin{Pf}
Suppose that $z$ is a central element and write
$$z=\sum_{\mu(M)<\lambda}c_M[M]+\sum_{\mu(M)=\lambda}c_M[M].$$
Let $N$ be any module and set $\nu:=\mu(N)$. Using Lemma \ref{ext} we know that
$$z[N]=\sum_{\mu(X)<\lambda\cup\nu}c_X'[X]+\sum_{\mu(M)=\lambda}c_MF_{MN}^{M\oplus N}[M\oplus N]$$
and similarly
$$[N]z=\sum_{\mu(X)<\lambda\cup\nu}c_X''[X]+\sum_{\mu(M)=\lambda}c_MF_{NM}^{M\oplus N}[M\oplus N].$$
Since
$$F_{MN}^{M\oplus N}=\frac{\left|\Aut(M\oplus N)\right|}{\left|\Aut M\right|\left|\Aut N\right|\left|\Hom(M,N)\right|},$$
we have that $F_{MN}^{M\oplus N}=F_{NM}^{M\oplus N}$ for all $N$ if and only if $\dim\Hom(M,N)=\dim\Hom(N,M)$ for all $N$. Applying the previous lemma, this is equivalent to $M=Z_\pi$ for some partition $\pi$, and so $\lambda=\pi^{\cup(l+1)}$.

Therefore there is a unique such module $M$ with $\mu(M)=\lambda$ appearing in the expression for $z$. By subtracting a suitable multiple of the central element $z_\pi=\prod_t z_{\pi_t}$ we obtain another central element $z'$, all of whose terms satisfy $\mu(M)\leq\lambda'<\lambda$. The proof now follows by induction on $\lambda$.
\end{Pf}

This completes the proof of Theorem \ref{ThmA}.

Fixing a simple $S$, the subalgebra of $\HA(\widetilde{\mathbb A}_{\,l})$ generated by the isomorphism classes $[S(ln)]$ for $n\geq0$ is isomorphic to $\Gamma$ via the homomorphism
$$\psi^{(S)}:\Gamma\to\HA(\widetilde{\mathbb A}_{\,l}), \quad \psi^{(S)}(P_\lambda(X;q^{-1}))=[S(l\lambda)],$$
where $S(l\lambda)=\bigoplus_i S(l\lambda_i)$.

\begin{Prop}\label{powersum}
For simples $S$ and $T$, we have that
$$\psi^{(S)}(p(n))-\psi^{(T)}(p(n))\in\CA(\widetilde{\mathbb A}_{\,l-1}).$$
\end{Prop}

I do not know a direct proof of this result, but we shall deduce it from the results of the following section on Hall algebras of affine quivers without oriented cycles.

\section{The affine quiver case}

Let $Q$ be an affine quiver without oriented cycles and $k$ a finite field with $q$ elements. Set $\HA$ and $\CA$ to be respectively the Hall algebra and the composition algebra of $kQ$. For a detailed description of the module category $\bmod kQ$ we refer the reader to \cite{DR}.

We recall that the tubes of the Auslander-Reiten quiver of the path algebra $kQ$ are indexed by the points of the scheme $\mathbb P^1_k$. Of these tubes, almost all are homogeneous. Write $dx$ for the degree of a point $x\in\mathbb P^1_k$. 

We fix a simple injective $I$ of defect $\langle\delta,\dimv I\rangle=1$. For each tube $x$ there exists a unique simple regular module $S_x$ mapping onto $I$. Set $R_x(n)$ to be the indecomposable regular module of dimension vector $ndx\delta$ into which $S_x$ embeds, and write $u_x(n)$ for the corresponding isomorphism class. Similarly we define $u_x(\lambda)$ to be the isomorphism class of $\bigoplus_i R_x(\lambda_i)$.

For the tube corresponding to $x$, let $\HA_x\subset\HA$ be the subalgebra generated by the isomorphism classes of indecomposables in this tube. We define an algebra homomorphism $\psi_x=\psi_x^{(I)}:\Gamma\to\HA_x$ as follows.

If the tube is homogeneous, then $\HA_x$ is isomorphic to the algebra $\HA_{dx}$ defined in Section \ref{HA} via $u_x(\lambda)\leftrightarrow u_{dx}(\lambda)$. Then $\psi_x$ is the isomorphism given by $\psi_x(P_\lambda(X;q^{-dx}))=q^{n(\lambda)dx}u_x(\lambda)$.

If the tube has rank $l\geq2$, then $dx=1$ and $\HA_x\cong\HA(\widetilde{\mathbb A}_{\,l-1})$. Inside $\HA_x$, though, we have a copy of $\HA_1$ given by $u_1(\lambda)\mapsto u_x(\lambda)$. Therefore we have the algebra monomorphism $\psi_x(P_\lambda(X;q^{-1}))=q^{n(\lambda)}u_x(\lambda)$.

\begin{Thm}\label{ThmB}
For $n\geq1$ and $P$ the indecomposable preprojective module of dimension vector $\dimv P=n\delta-\dimv I$, we have
$$[I][P]=q^{n-1}[P][I]+\frac{1}{q-1}\mathfrak c(n).$$
The $\mathfrak c(n)$ are given by
$$\mathfrak c(n)=\sum \left|\Aut X\right|[X]$$
where the sum is taken over all modules $X=\bigoplus_i R_{x_i}(m_i)$ such that $\dimv X=n\delta$ and there is at most one summand from each tube.
\end{Thm}

This is a generalisation of Theorem 4.3 in \cite{Szanto}, which deals with the Kronecker quiver.

We shall divide the proof into a sequence of easy lemmas.

\begin{Lem}
Let $X=R_{x_1}(m_1)\oplus\cdots\oplus R_{x_r}(m_r)$ with the $x_i$ distinct points of $\mathbb P^1_k$ and $\sum_im_idx_i=n$. Then the Hall number $F_{IP}^X\neq0$.
\end{Lem}

\begin{Pf}
Suppose first that $X=R_x(m)$. Then there is a unique exact sequence
$$0\to S_x\to X\to T\to 0.$$
Since $\dim\Hom(X,I)=n$ whereas $\dim\Hom(T,I)=n-1$, we can lift the canonical epimorphism $S_x\to I$ to an epimorphism $g:X\to I$ such that $g$ does not factor through any proper regular factor of $X$.

The kernel $K$ must be indecomposable preprojective, hence isomorphic to $P$. For, $K$ cannot contain any preinjective direct summand. Moreover, since the defect of $K$ is $-1$, there is precisely one preprojective direct summand. Suppose that there exists a regular direct summand $R'$ of $K$. Then we have the commutative diagram
$$\begin{CD}
0 @>>> R' @>>> X @>>> R'' @>>> 0\\
&& @VVV @| @VVV\\
0 @>>> K @>>> X @>>> I @>>> 0
\end{CD}$$
It follows that the map $X\to I$ factors through $R''$, a contradiction.

Now suppose that $X=\bigoplus_i R_{x_i}(m_i)$. For each summand we have an epimorphism $g_i:R_{x_i}(m_i)\to I$ with indecomposable preprojective kernel, so consider the map $(g_i):X\to I$. Again the kernel $K$ contains no preinjective direct summand and has defect $-1$, so there is precisely one preprojective summand.

Suppose that the kernel decomposes as $K=K'\oplus R'$ with $R'$ indecomposable regular. Then we can write $X=X'\oplus R_x(m)$ with $\Hom(R',X')=0$. It now follows that the composition $R'\to R_x(m)\to I$ is zero, and so the image of $R'$ (which is regular) lies in the kernel of $R_x(m)\to I$ (which is preprojective), a contradiction.
\end{Pf}

We shall need the following observation, which follows from the $k$-vector space structure on $\Ext^1_R(M,N)$.
\begin{Lem}
Let $k$ be a field and $R$ a $k$-algebra. For $R$-modules $M$ and $N$, the subgroup $k^*$ of $R$ acts freely on $\Ext^1_R(M,N)-\{0\}$.
\end{Lem}

We recall the following formula of Riedtmann \cite{Riedt}
$$F_{MN}^X=\varepsilon_{MN}^X\frac{\left|\Aut X\right|}{\left|\Aut M\right|\left|\Aut N\right|\left|\Hom(M,N)\right|}$$
expressing the Hall number $F_{MN}^X$ in terms of the number $\varepsilon_{MN}^X=\left|\Ext^1(M,N)_X\right|$ of extension classes with middle term $X$.

In our situation we have that $F_{IP}^X=\varepsilon_{IP}^X\frac{\left|\Aut X\right|}{(q-1)^2}$.

Now, for $X=\bigoplus_i R_{x_i}(m_i)$, the Hall number $F_{IP}^X$ is non-zero by the first lemma, and hence the number $\varepsilon_{IP}^X$ has size at least $q-1$. We also know that
$$\dim\Ext^1(I,P)=-\langle\dimv I,\dimv P\rangle=-\langle\dimv I,n\delta-\dimv I\rangle=n+1.$$

\begin{Lem}
We have the identity
$$\sum_{\substack{(x_1,m_1),\ldots,(x_r,m_r)\\x_i\ \mathrm{distinct}\\\sum_im_idx_i=n}}\!\!\!\!\!\!\!\!\!\! 1 \quad =\frac{q^{n+1}-1}{q-1}.$$
\end{Lem}

\begin{Pf}
Consider the generating function for these numbers, namely
$$N(t)=\prod_x(1+t^{dx}+t^{2dx}+\cdots)=\prod_x\frac{1}{1-t^{dx}}=\prod_{d\geq1}\Big(\frac{1}{1-t^d}\Big)^{\phi(d)},$$
where $\phi(d)$ is the number of elements of degree $d$ in $\mathbb P^1_k$. That is, $\phi(1)=q+1$ and $\phi(d)$ for $d\geq 2$ is the number of monic polynomials of degree $d$ over $k$.

Then
\begin{align*}
\log N(t) &= -\sum_{d\geq1}\phi(d)\log(1-t^d) = \sum_{d\geq1}\phi(d)\sum_{r\geq1}\frac{1}{r}t^{rd} = \sum_{n\geq1}\frac{1}{n}\Big(\sum_{d|n}d\phi(d)\Big)t^n\\
& = \sum_{n\geq1}\frac{1}{n}(q^n+1)t^n = -\log(1-tq)-\log(1-t).
\end{align*}
Therefore $N(t)=(1-t)^{-1}(1-tq)^{-1}=\sum_{n\geq0}(q^n+\cdots+q+1)t^n$.
\end{Pf}

It follows that the number of choices for $X=\bigoplus_i R_{x_i}(m_i)$ of dimension vector $n\delta$ with the $x_i$ distinct is $(q^{n+1}-1)/(q-1)$ and thus we have accounted for at least $q^{n+1}-1$ extension classes. Once we have included the trivial class, which satisfies $\varepsilon_{IP}^{P\oplus I}=1$, we see that we are done.

In summary, we know that $\varepsilon_{IP}^X=q-1$ and $F_{IP}^X=\left|\Aut X\right|/(q-1)$ for each $X=\bigoplus_i R_{x_i}(m_i)$ of dimension vector $n\delta$ such that the $x_i$ are distinct, and that $F_{IP}^{P\oplus I}=\left|\Hom(P,I)\right|=q^{n-1}$. Moreover, these are all the terms which occur in the product $[I][P]$.

The proof of Theorem \ref{ThmB} is now complete.

Let $x_1,\ldots,x_r$ represent the non-homogeneous tubes and set $\CA_{x_i}\subset\HA_{x_i}$ to be the corresponding composition algebra. Therefore $\CA_{x_i}$ is the subalgebra generated by the isomorphism classes of the regular simples in the $x_i$-th tube.

\begin{Prop}
The composition algebra $\CA$ of $kQ$ contains the subalgebras
$$\bQ[\mathfrak p(1),\mathfrak p(2),\ldots] \quad\textrm{and each}\quad \CA_{x_1},\ldots,\CA_{x_r},$$
where the $\mathfrak p(n)$ are given by the formula
$$\mathfrak p(n)=\sum_{mdx=n}\frac{1}{m}\psi_x(p(m))=\sum_{mdx=n}\frac{1}{m}\sum_{\lambda\vdash m}(1-q^{dx})\cdots(1-q^{(\ell(\lambda)-1)dx})u_x(\lambda).$$
\end{Prop}

\begin{Pf}
We have just shown that the elements $\mathfrak c(n)$ all lie in the composition algebra. Let $C(t)=1+\sum_{n\geq1}\mathfrak c(n)t^n$ be their generating function. Then
\begin{align*}
C(t) &= \prod_x\Big(1+\sum_{n\geq1}(1-q^{-dx})q^{ndx}u_x(n)t^{ndx}\Big)\\
&= \prod_x A_{dx}(t) = \prod_x \psi_x\big(C(q^{-dx},(qt)^{dx})\big).
\end{align*}
In particular,
\begin{align*}
\diff\log C(t) &= \sum_x\psi_x\Big(\diff\log C(q^{-dx},(qt)^{dx})\Big)\\
&= \sum_x dx\,t^{dx-1}\sum_{n\geq1} (q^{ndx}-1)\psi_x(p(n))t^{(n-1)dx}\\
&= \sum_x dx\sum_{n\geq1} (q^{ndx}-1)\psi_x(p(n))t^{ndx-1}\\
&= \sum_{n\geq1} n(q^n-1)\mathfrak p(n)t^{n-1},
\end{align*}
and so each $\mathfrak p(n)$ also lies in $\CA$.

Now, any algebraic relation satisfied by the $\mathfrak p(n)$ would also be satisfied by the $\psi_x(p(n))$ for each $x$ of degree 1, and hence for the $p(n)$ themselves, which we know are algebraically independent. Hence the $\mathfrak p(n)$ are algebraically independent and
$$\bQ[\mathfrak c(1),\mathfrak c(2),\ldots]=\bQ[\mathfrak p(1),\mathfrak p(2),\ldots].$$

For the second part, we note that each regular simple from a non-homogeneous tube has dimension vector a real root less than $\delta$. Since there are only finitely many isomorphism classes of indecomposables of dimension vector less than $\delta$, and each is completely determined by its dimension vector, it follows (analogously to the arguments in the Dynkin quiver case, c.f. \cite{Ringel3}) that the isomorphism class of each such regular simple lies in $\CA$. It follows that each $\CA_x$ for $x$ a non-homogeneous tube is a subalgebra of $\CA$.
\end{Pf}

We note that the element $\mathfrak p(n)$ has the nice property that each isomorphism class which occurs in its expression lies in a single tube.

For the preprojective component of the Auslander-Reiten quiver of $Q$, there is a natural partial ordering whereby
$$\dim\Hom(P,P')\neq0 \quad\textrm{implies}\quad [P]\preceq[P'].$$
We complete this to a total ordering $\leq$. Similarly for the preinjective component.

We can also choose a PBW basis for each composition algebra $\CA_x$ with $x$ non-homogeneous. For example, such a basis is provided by Theorem 9.3 of \cite{DD}.

We can now state our main theorem of this section. This is a generalisation of Theorem 6.1 in \cite{Szanto}, and a refinement of the Main Theorem in \cite{Zhang}.

\begin{Thm}\label{ThmC}
The composition algebra $\CA$ has a PBW basis consisting of elements of the form $PRC_1\cdots C_rI$, where
\begin{enumerate}
\item $P=[P_1]\cdots[P_r]$ with $P_1\leq\cdots\leq P_r$ indecomposable preprojectives;
\item $R=\mathfrak p(a_1)\cdots\mathfrak p(a_s)$ with $a_1\leq\cdots\leq a_s$;
\item $C_i$ is a PBW-basis element for the composition algebra $\CA_{x_i}$;
\item $I=[I_1]\cdots[I_t]$ with $I_1\leq\cdots\leq I_t$ indecomposable preinjectives.
\end{enumerate}
In particular, we have the decomposition
$$\CA=\mathcal P\cdot\bQ[\mathfrak p(1),\ldots]\cdot\CA_{x_1}\cdots\CA_{x_r}\cdot\mathcal I.$$
\end{Thm}

\begin{Pf}
We know that $\mathcal P$ and $\mathcal I$ are subalgebras of the composition algebra from the results in Section \ref{RF}. Furthermore, the previous proposition tells us that $\bQ[\mathfrak p(1),\ldots]$ as well as each $\CA_{x_i}$ are all subalgebras of $\CA$. The right hand side is therefore a graded vector subspace of the composition algebra and so it is enough to show that the dimensions of each graded part coincide.

Let $\mathfrak g$ be the affine Kac-Moody Lie algebra corresponding to the underlying graph of $Q$. The positive roots of $\mathfrak g$ are of two types \cite{Kac3}: the real roots (having multiplicity one) and the multiples of $\delta$ (having multiplicity $N$, the number of vertices of $Q$).

We also recall that the dimension vectors of the indecomposable modules for $Q$ are precisely the positive roots of $\mathfrak g$. (This is proved in \cite{DR} for the affine case; see also \cite{Kac1, Kac2, KR} for the general result.)

It follows from the PBW basis for $\QG^+(\mathfrak g)$ that the dimension of the graded part $\CA_\alpha$ is precisely the number of ways of expressing $\alpha$ as a sum of positive roots (with multiplicity). So for example, $\dim\CA_\delta$ equals $N$ plus the number of ways of expressing $\delta$ as a sum of positive real roots.

Similarly, we can apply this same result to each $\CA_x$ for $x$ a non-homogeneous tube of rank $l$. The corresponding Lie algebra is then $\mathfrak g_{\,x}:=\widehat{\mathfrak{sl}}_{\,l}$. The roots of $\mathfrak g_{\,x}$ correspond to the dimension vectors of the indecomposables in the $x$-th tube, with the simple roots being sent to the dimension vectors of the corresponding regular simples. In particular, $m\delta$ has multiplicity $l-1$.

Moreover, Theorem 4.1 in \cite{DR} states that the sum $\sum_i(l_i-1)$ of the ranks of the non-homogeneous tubes minus 1 equals $N-1$, the number of vertices of $Q$ minus 1.

Now, for an arbitrary preprojective module $P$ we can write $P=P_1^{a_1}\oplus\cdots\oplus P_n^{a_n}$ such that $[P_1]<\cdots<[P_n]$ are indecomposable. Then
$$[P]=\gamma[P_1]^{a_1}\cdots[P_n]^{a_n}, \qquad\textrm{for some } \gamma\neq0.$$
It follows that $\dim\mathcal P_\alpha$ equals the number of ways of expressing $\alpha$ as a sum of dimension vectors of indecomposable preprojective modules. An analogous result holds for the subalgebra $\mathcal I$ of $\CA$ consisting of all preinjective modules.

Finally, each $\mathfrak p(n)$ is homogeneous of degree $n\delta$.

We deduce that the dimensions of the homogeneous parts of degree $\alpha$ of the subspace $\mathcal P\cdot\bQ[\mathfrak p(1),\ldots]\cdot\CA_{x_1}\cdots\CA_{x_r}\cdot\mathcal I$ and of $\CA$ coincide. This completes the proof of the theorem.
\end{Pf}

We now make the following observation.

\begin{Cor}
Let $I$ and $J$ be simple injectives of defect 1. We can use these to define morphisms $\psi_x^{(I)}$ and $\psi_x^{(J)}$, and hence obtain elements $\mathfrak p^{(I)}(n)$ and $\mathfrak p^{(J)}(n)$. These are then related by
$$\mathfrak p^{(I)}(n)-\mathfrak p^{(J)}(n)\in\CA_{x_1}\oplus\cdots\oplus\CA_{x_r}.$$
In particular, we deduce Proposition \ref{powersum} of the previous section.
\end{Cor}

\begin{Pf}
We know that $\mathfrak p^{(I)}(n)=\sum_{mdx=n}\frac{1}{m}\psi_x^{(I)}(p(m))$, where $\psi_x^{(I)}$ is defined with respect to the regular simple $S_x\twoheadrightarrow I$. Similarly $\psi_x^{(J)}$ is defined using $T_x\twoheadrightarrow J$.

For $x$ a homogeneous tube, there exists a unique regular simple and hence the two maps $\psi_x^{(I)}$ and $\psi_x^{(J)}$ coincide. Therefore the difference $\mathfrak p^{(I)}(n)-\mathfrak p^{(J)}(n)$ is a sum of elements lying in the Hall algebras $\HA_{x_i}$ of the non-homogeneous tubes, and is also an element of the composition algebra $\CA$. It follows that this difference is a sum of elements lying in the composition algebras $\CA_{x_i}$ of the non-homogeneous tubes. That is,
$$\psi_x^{(I)}(p(n))-\psi_x^{(J)}(p(n))\in\CA_x$$
for each non-homogeneous $x$.

Now consider the quiver of type $\widetilde{\mathbb A}_{\,l+2}$ with orientation
$$
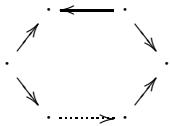
\xymatrix@R=10pt@C=7pt{&\cdot&&\cdot\ar[ll]\ar[dr]\\\cdot\ar[ur]\ar[dr]&&&&\cdot\\&\cdot\ar@{.>}[rr]&&\cdot\ar[ur]}$$
This has two simple indecomposables $I$ and $J$ of defect 1, and there are two non-homogeneous tubes, of ranks 2 and $l+1$. Moreover, without loss of generality we may suppose that we have a regular simple $S$ in the tube of rank $l$ such that ${}^\tau S\twoheadrightarrow I$ and ${}^{\tau^2} S\twoheadrightarrow J$.

It follows that, in the notation of Proposition \ref{powersum}, the difference $\psi^{(S)}(p(n))-\psi^{(T)}(p(n))$ for $T={}^\tau$ lies in the composition algebra for a tube of rank $l+1$, or equivalently in $\CA(\widetilde{\mathbb A}_{\,l})$, and thus the same holds for all regular simples $S$ and $T$. This proves Proposition \ref{powersum}.
\end{Pf}


\begin{thebibliography}{99}
\bibitem{DD} \textsc{B.~Deng} and \textsc{J.~Du}, `Monomial bases for quantum affine $\mathfrak{sl}_{\,n}$' (preprint) math.RA/0307258 .
\bibitem{DR} \textsc{V.~Dlab} and \textsc{C.M.~Ringel}, `Indecomposable representations of groups and algebras', \textit{Mem. Amer. Math. Soc.} 6 (1976).
\bibitem{Green} \textsc{J.M.~Green}, `Hall algebras, hereditary algebras and quantum groups', \textit{Invent. Math.} 120 (1995) 361--377.
\bibitem{Hall} \textsc{P.~Hall}, `The algebra of partitions', \textit{Proc. 4th Canadian Math. Congress} (1959), 147--159.
\bibitem{Kac1} \textsc{V.G.~Kac}, `Infinite root systems, representations of graphs and invariant theory', \textit{Invent. Math.} 56 (1980) 57--92.
\bibitem{Kac2} \textsc{V.G.~KAc}, `Root systems, representations of quivers and invariant theory', \textit{Invariant theory} (ed. F.~Gherardelli), Lecture Notes in Mathematics 996 (Springer, Berlin, 1983) 74--108.
\bibitem{Kac3} \textsc{V.G.~Kac}, \textit{Infinite dimensional Lie algebras} (3rd ed.) (Cambridge University Press, Cambridge, 1990).
\bibitem{KR} \textsc{H.~Kraft} and \textsc{C.~Riedtmann}, `Geometry of representations of quivers', \textit{Representations of algebras} (ed. P.~webb), London Math. Soc. Lecture Notes Series 116 (Cambridge University Press, Cambridge, 1986) 109--145.
\bibitem{Lusztig} \textsc{G.~Lusztig}, \textit{Introduction to quantum groups} (Birkh\"auser, Boston, 1993).
\bibitem{Lusztig1} \textsc{G.~Lusztig}, `Canonical bases arising from quantized enveloping algebras', \textit{J. Amer. Math. Soc.} 3 (1990) 447--498.
\bibitem{Macd} \textsc{I.G.~Macdonald}, \textit{Symmetric functions and Hall polynomials} (2nd ed.), Oxford Mathematical Monographs (Clarendon Press, Oxford, 1995).
\bibitem{Reineke} \textsc{M.~Reineke}, `Generic extensions and multiplicative bases of quantum groups at $q=0$', \textit{Represnt. Theory} 5 (2001) 147--163.
\bibitem{Riedt} \textsc{C.~Riedtmann}, `Lie algebras generated by indecomposables', \textit{J. Algebra} 170 (1994) 526--546.
\bibitem{Ringel1} \textsc{C.M.~Ringel}, `Hall algebras and quantum groups', \textit{Invent. Math.} 101 (1990) 583--592.
\bibitem{Ringel2} \textsc{C.M.~Ringel}, `The composition algebra of a cyclic quiver. Towards an explicit description of the quantum group of type $\tilde A_n$', \textit{Proc. London Math. Soc.} (3) 66 (1993) 507--537.
\bibitem{Ringel3} \textsc{C.M.~Ringel}, `The Hall algebra approach to quantum groups' (E.L.A.M. Lectures 1993) \textit{Aportaciones Mat. Comun.} 15 (1995) 85--114. 
\bibitem{Ringel4} \textsc{C.M.~Ringel}, `Green's theorem on Hall algebras', \textit{Representation theory of algebras}, CMS Conference Proceedings 19 (Amer. Math. Soc., Providence, 1996) 185--245.
\bibitem{Schiff} \textsc{O.~Schiffmann}, `The Hall algebra of a cyclic quiver and canonical bases of Fock spaces', \textit{Intern. Math. Res. Notices} 8 (2002) 413--440.
\bibitem{Stein} \textsc{E.~Steinitz}, `Zur Theorie der Abel´schen Gruppen', \textit{Jahresbericht der DMV} 9 (1901) 80--5.
\bibitem{Szanto} \textsc{C.~Sz\'ant\'o}, `Hall numbers and the composition algebra of the Kronecker algebra' (preprint) http://math.ubbcluj.ro/\~{}szanto/publ/KronHall1.pdf .
\bibitem{Zhang} \textsc{P.~Zhang}, `Composition algebras of affine type', \textit{J. Algebra} 206 (1998) 505--540. 
\end{thebibliography}
\end{document}